\documentclass[a4paper, 11pt]{amsart}
\usepackage{amsmath,amssymb,amsthm}
\usepackage{a4wide}
\usepackage{MnSymbol}

\newtheorem{theorem}{Theorem}

\newtheorem{corollary}[theorem]{Corollary}
\theoremstyle{definition}

\newtheorem{example}[theorem]{Example}
\newtheorem{lemma}{Lemma}

\newtheorem{remark}[theorem]{Remark}
\newtheorem{proposition}[theorem]{Proposition}

\newcommand{\pd}[2]{\frac{\partial#1}{\partial#2}}
\newcommand{\vp}[1]{\frac{\partial}{\partial#1}}

\newcommand{\X}{\mathfrak{X}}

\renewcommand{\d}{\text{d}}
\newcommand{\Sp}{\text{Sp}}
\newcommand{\R}{\mathbb{R}}

\begin{document}
\title{Integration of PDEs by differential geometric means}
\author{Naghmana Tehseen and Geoff Prince}
\address{ Department of Mathematics and Statistics, La Trobe University, Victoria, 3086, Australia}
\email{ntehseen@students.latrobe.edu.au,\ g.prince@latrobe.edu.au}
\date{October 2012}
\noindent \begin{abstract}
We use Vessiot theory and exterior calculus to solve partial differential equations(PDEs) of the type $u_{yy}=F(x,y,u,u_x,u_y,u_{xx},u_{xy})$ and associated evolution equations.  These equations are represented by the Vessiot distribution of vector fields. We develop and apply an algorithm to find the largest integrable sub-distributions and hence solutions of the PDEs.  We then apply the integrating factor technique \cite{sherr} to integrate this integrable Vessiot sub-distribution. The method is successfully applied to a large class of linear and non-linear PDEs.
 \end{abstract}
\maketitle
\section{Introduction}
\noindent The geometric study of differential equations aims to describe the local and global structure of the
solutions as integral submanifolds of an equation manifold. This study of PDEs leads to the concept of {\em Vessiot distributions}, which provides the convenient formal framework to investigate PDEs on an appropriate jet bundle.

In \cite{sherr} the authors used exterior calculus and the Frobenius theory of foliations to show that integrating factors can be used to solve systems of higher order ordinary differential equations, and they considerably extended the class of group actions which could be thought of as symmetries of the system for the purpose of reduction to the quadratures. Moreover their methods avoided cumbersome changes of coordinates, allowing solutions and first integrals to be represented in the original, physical coordinates of the problem. Later this technique is extended to the first order PDEs in \cite{Mbarco01,Mbarco02,MbarcoGP01,MBGP01}. In the present paper we will extend the integrating factor technique to second order PDEs  by using their results  to integrate the (integrable sub-distributions of the) vector field description of the PDEs.

The formulation of a (finite) system of differential equations as a distribution of vector fields on some appropriate jet space is due to Vessiot \cite{Vessiot1924,Vessiot1939,Vessiot1942} and called a Vessiot distribution. The maximal integrable sub-distributions of the  Vessiot distribution represent every local solution of the PDE system. Computing the Vessiot distribution is routine but finding the maximal integrable sub-distributions is not. In fact there is no known algorithm for computing them in the generic case.

Solving the PDE system is equivalent to integrating a Pfaffian system or its dual vector field system.
The Vessiot formulation of PDEs is dual to the more popular {\em exterior differential system (EDS)} formulation due to Cartan \cite{Cartan1945} and elaborated in \cite{Ivey,bryant}. In \cite{Fackerell} the author pointed out that the former is in many ways computationally  simpler because vector fields are simpler to implement than the full exterior calculus of $p-$forms. Recent developments in the Vessiot theory can be found in \cite{Fesser09,Stormark,Vassiliou01}.

The purpose of this paper is to provide a systematic approach to solving second order linear and non-linear PDEs in the presence of symmetries. We use a general class of symmetries known as \emph{solvable structures} which are not necessarily of point type. To generate the solvable structures for each integrable sub-distributions we use the symmetry determination software package DIMSYM \cite{sherrdimsym} operating as a REDUCE \cite{Reduce}  overlay. We also use the exterior calculus package EXCALC \cite{Sexcalc}.

The paper is organised as follows. In Section $2$ we introduce the key concepts of the Vessiot distribution. In Section $3$ we give an explanation of how a reduction of order is possible using symmetries and present the main result of \cite{sherr} which we will use to integrate the Frobenius integrable distribution of vector fields. In Section $4$ we describe, in detail, a three steps algorithm to find the largest integrable sub-distributions which satisfy the appropriate independence/projectability  condition. To integrate such distributions we use the results of \cite{sherr}. We present some examples to explain the algorithm. In Section $5$ we produce an EDS method to solve the differential conditions in step $2$ of the algorithm. We also demonstrate the construction of a group invariant solution of our second order PDEs. We will explain the methods with some concrete examples.

\section{Basic concepts and notations}
\noindent Consider a system of partial differential equations  \cite{Stormark} of $m$ independent variables $x^i$ and $n$ dependent variables  $u^j$,
\begin{align}F^a(x^i,u^j,u^j_{i_1},u^j_{i{_1}i{_2}}\ldots u^j_{i_1\ldots i_k})=0. \label{eqn:g.pde}\end{align} Denote by $X$ the space of independent variables and by $U$ the space of dependent variables. The subscripts $1\leq i_1\leq\ldots \leq i_k\leq n$ are used to specify the partial derivatives of $u^j$, where $k$ is the maximum order of the system. The reader may consult \cite{Stormark} for more details.

\noindent Consider the trivial bundle $\pi:X\times U\rightarrow X$ and let $f$ and $g$ be two smooth sections of $\pi$. We say that $f$ and $g$  are equivalent to order $k$ at $x$ if and only if  \[\frac{\partial^{q_1+ \ldots +q_m} f}{(\partial x^1)^{q_1}\ldots(\partial x^m)^{q_m}}(x)=\frac{\partial^{q_1+\ldots+q_m} g}{(\partial x^1)^{q_1}\ldots(\partial x^m)^{q_m}}(x)\] for all $m-$tuples $(q_1\ldots q_m)$ with $q_1+\ldots+q_m\leq k.$
The equivalence class of a smooth section $f$ at a point $x$ is the $k-$jet of $f$ at $x$  and is denoted by $j^k_x f.$ The set of all $k-$jets of smooth sections of $\pi$ constitutes the bundle of $k-$jets of maps $X\rightarrow U$ and denoted by $J^k(X,U).$ The zeroth order jet bundle, $J^0(X,U),$ is identified with $X\times U.$ The $k-$graph of smooth section $f$ of $\pi$ is the map $j^k f: X\rightarrow J^k(X,U)$ defined by $x\rightarrow j^k_x f.$ Thus, the image of $k-$graph of a section is a $m-$dimensional immersed submanifold of $J^k(X,U).$

The geometry of jet bundles is to a large extent determined by their \emph{contact structure.} Each jet bundle $J^k(X,U)$ comes equipped with a module of differential $1-$forms spanned by  \[\theta^j_I:=d u^j_I-\sum_{i=1}^m u_{I, i}^j d x^i,\] where $I$ is a multi-index of order less than or equal to $k-1$ and is denoted by $\Omega^k(X,U).$ The contact system $\Omega^k(X,U)$ on $J^k$ is also an example of an exterior differential system generating an \emph{ideal} in the ring of differential forms on a differentiable manifold. This particular differential system generated by $1-$forms is also called a {\em Pfaffian system}.

Another way to view this is that there may be a submanifold $N \subset M$ where all the forms in the contact co-distribution are zero when pulled back to $N$. That is, if $i:X \rightarrow J^k(X,U)$ is any immersion whose pull-back annihilates the contact system and satisfies the {\em transverse condition} or {\em independence condition} $i^\ast(dx^1\wedge\ldots \wedge dx^n)\neq 0,$ then $i(X)$ is the image of some $k-$jet. The image of the $k-$jet of any function $f:X\rightarrow U$ is an {\em integral submanifold} of the $k$th order contact system.

Now we can see how these constructions are related to differential equations. The study of second and higher order differential equations on manifolds needs to consider jet bundles on some underlying space. In view of the definition of a jet bundle, the $k-$jet bundle $J^k(X,U)$ is the appropriate framework for dealing with the systems of PDEs of a higher order of $m$ dependent variables and $n$ independent variables. We will be interested in the second order partial differential equations which define embedded submanifolds of $J^k(X,U).$

\noindent A solution of \eqref{eqn:g.pde}  is a smooth function $f:X\rightarrow U$ whose $k-$graph $J^k_j$ defines a $m-$dimensional submanifold $S$ of a smooth embedded submanifold of $J^k(X,U)$.

\section{Ideals, Symmetries and reduction of order}
\noindent We will briefly explain how the integrating factor technique in \cite{sherr} allows a reduction of order via symmetries for higher order ordinary differential equations.
We begin with some basics, for more details see \cite{MbarcoGP01,MBGP01,cramp,bryant,sherr}.

\noindent Let $M$ be some smooth manifold of dimension $n.$ For simplicity, we suppose that all our objects are smooth on $M$. The algebra of differential forms on $M,~\bigwedge(M),$ is a graded algebra. An \textit{ideal} $I$ in $\bigwedge(M)$ is an additive subgroup of $\bigwedge(M)$ that is closed under the wedge product $(\beta^i\in \textit{I}$ implies $\beta^i\wedge \alpha \in \textit{I},\forall \alpha \in \bigwedge(M)).$  An ideal \emph{I} is a {\em differential ideal} if the exterior derivative of every member of \emph{I} is also in \emph{I}. A vector field $V$ is said to be a {\em symmetry} of an ideal if $\mathcal{L}_V \emph{I}\subset \emph{I}.$

 The {\em kernel} of a differential form $\Omega$ is the submodule of vector fields annihilating $\Omega$, $\text{ker}\ \Omega:=\{Y\in\X(M):Y\righthalfcup \Omega=0\}$.
A differential $p-$form $\Omega$ on $M$ is {\em simple or decomposable} if it is the wedge product of $p$ $1-$forms, so $\Omega$ is locally simple if $\text{ker}\ \Omega $ is everywhere $n-p$ dimensional, that is, of maximal dimension. A {\em constraint} $1-$form $\theta$ for differential form $\Omega$ is a $1-$form satisfying $\theta\wedge\Omega=0,$ which implies $Y\righthalfcup \theta=0,\ \forall ~Y\ \in \text{ker}~\Omega.$ An $m$-dimensional distribution $D$ on $M$ is an assignment of an $m$-dimensional subspace $D_x$ of $T_xM$ to each point $x$ of $M$. The corresponding co-distribution $D^\perp$ consists of the annihilators $D_x^\perp \subset T_x^\ast M$. As a submodule of $\bigwedge^1 (M)$ the elements of $D^\perp$ are called constraint forms for $D$. A {\em characterising} form for an $m-$dimensional distribution $D$ is a form $\Omega$ on $M$ of degree $(n-m)$ which is the exterior product of $(n-m)$ constraint forms.

\noindent A distribution $D$  has (maximal) integral submanifolds if and only if\begin{align}\label{eqn:FI}
&d\theta\wedge \Omega=0, \ \ \ \ \forall\ \  \theta \in D^{\perp}\ \ (\Omega \ \text{is any characterising form for D})\\ \text {or}\nonumber\\
&[Y,Z]\in D,\ \ \ \ \ \forall Y, Z \in D.\nonumber\end{align} In this case $D$ is said to be \emph{involutive or Frobenius integrable} and any characterising form $\Omega$ is said to be Frobenius integrable.
 $\Omega$ is Frobenius integrable if its kernel is Frobenius integrable and of maximal dimension everywhere. A function $f$ on $M$ is called an {\em integrating factor} for the Frobenius integrable $1-$form $\theta$ if $d(f\theta)=0.$ A vector field $V\in \X(M)$ is a symmetry of a vector field distribution $D\subset \X(M) $ if $\mathcal{L}_V D\subset D$ and $V$ is a symmetry of differential form $\Omega$ if $\mathcal{L}_V \Omega=\lambda \Omega \ \ \text{for some smooth function}\ \ \lambda\ \text{on}\ M. $ If $V$ is in $D$ it is said to be {\em trivial}. For any vector field distribution $D,$  we say that a collection of $k$ linearly independent vector fields $X_1,\ldots,X_k\in\X(M)$ forms a {\em solvable structure} for $D$ if
  \begin{align*}
  \pounds_{X_k}D~&\subset~ D,\\
  \pounds_{X_{k-1}}(\Sp\{X_k\}\bigoplus D)~&\subset~\Sp\{X_k\}\bigoplus D,\\
  \vdots\nonumber\\
\pounds_{X_1}(\Sp\{X_2,\ldots,X_k\}\bigoplus D)~&\subset~\Sp\{X_2,\ldots,X_k\}\bigoplus D.
\end{align*}
DIMSYM can be used to locate solvable structures for distributions. The following results will be useful.

\begin{theorem}\cite{sherr}\label{thm:sher2.1}
Let $\theta$ be Frobenius integrable $1-$form which is nowhere zero on some open subset
$U$ of a manifold M. Then $\mu$ is an integrating factor for $\theta$
if and only if $\mu=(X\righthalfcup\theta)^{-1}$ for some symmetry vector field $X$
on $U$ with $X\righthalfcup\theta\neq 0$ on $U.$
\end{theorem}

\begin{proposition}\label{prop:d(a^b)=k^a^b}
Let $\theta^1,\theta^2,~\ldots,~\theta^k \in \bigwedge^1(M)$  then $\Sp\{\theta^1,\theta^2,~\ldots,~\theta^k\}$ is Frobenius integrable, \[d\theta^a\wedge\theta^1\wedge\ldots\wedge\theta^k =0,~a=1,\ldots,k,\]
if and only if $d(\theta^1\wedge\ldots\wedge\theta^k)= \lambda\wedge\theta^1\wedge\ldots\wedge\theta^k,\ \ \ \text{for some 1-form}\ \  \lambda.$
\begin{proof}
If $\Sp\{\theta^1,\theta^2,~\ldots,~\theta^k\}$ is Frobenius integrable, then $d\theta^a\wedge\theta^1\wedge\ldots\wedge\theta^k=0,$
and from the definition of exterior derivative $d(\theta^1\wedge\ldots\wedge\theta^k)=d\theta^1\wedge\theta^2\wedge\ldots\wedge\theta^k+\ldots+(-1)^k \theta^1\wedge\ldots\wedge d\theta^k$.  It is obvious from the Frobenius conditions that $(\theta^1\wedge\ldots\wedge\theta^k)$ divides $d(\theta^1\wedge\ldots\wedge\theta^k)$. Hence, $d(\theta^1\wedge\ldots\wedge\theta^k)= \lambda \wedge\theta^1\wedge\ldots\wedge\theta^k.$

Conversely, assume that $d(\theta^1\wedge\ldots\wedge\theta^k)= \lambda \wedge\theta^1\wedge\ldots\wedge\theta^k$ then, by using the definition of exterior derivative,
$$d\theta^1\wedge\theta^2\wedge\ldots\wedge\theta^k+\ldots+(-1)^k\theta^1\wedge\ldots\wedge d\theta^k= \lambda\wedge \theta^1\wedge\ldots\wedge\theta^k.$$ By taking the wedge product with $\theta^a,$ we have
\[d\theta^a\wedge\theta^1\wedge\ldots\wedge\theta^k=0.\] So,
$\Sp\{\theta^1,~\ldots,~\theta^k\}$ is Frobenius integrable.
\end{proof}
\end{proposition}

We now use Proposition \ref{prop:d(a^b)=k^a^b} to give a generalisation of Theorem~\ref{thm:sher2.1} (compare with \cite{MBGP01}).

\begin{theorem}\label{thm2}
Let $\Omega$ be a k-form on a manifold M, and
$\Sp(\{X_1\ldots X_k\})$ be a k-dimensional distribution on open
$U\subseteq M$ satisfying $X_i\righthalfcup \Omega\neq 0$ everywhere
on $U$. Further suppose that $\Sp(\{X_1\ldots X_k\}\cup
~\text{ker}~\Omega)$ is Frobenius integrable for some $j<k$ and that $X_i$ is
a symmetry of $\Sp(\{X_{i+1}\ldots X_k\}\cup ~\text{ker}~ \Omega) ~\text{for}~
i=1\ldots j.$

Put $\alpha^i:=X_1\righthalfcup\ldots \righthalfcup\bar{X}_i\righthalfcup\ldots \righthalfcup X_k
\righthalfcup \Omega$, where $\bar{X}_i$ indicates that this argument is missing, and put
$\omega^i:=\displaystyle{\frac{\alpha^i}{X_i\righthalfcup \alpha^i}} \ \text{for}\
i=1\ldots k$, so that $\{\omega^1\ldots \omega^k\}$ is dual to
$\{X_1\ldots X_k\}$.
Then\begin{enumerate}
\item[(i)]$\omega^1\wedge \omega^2\wedge
\ldots \wedge\omega^k=\displaystyle{\frac{\Omega}{\Omega(X_1,X_2\ldots X_k)}}$
\item[(ii)] $\omega^1\wedge\omega^2\wedge\ldots \wedge\omega^k$ is closed.
\end{enumerate}\begin{proof}
 To prove the first part, observe that
 \[\omega^1\wedge\omega^2=-~\frac{(X_1\righthalfcup\Omega)\wedge(X_2\righthalfcup\Omega)}{(\Omega(X_2,X_1))^2}.\]
 We know that $(X_1\righthalfcup\Omega)\wedge(X_2\righthalfcup\Omega)
 =X_1\righthalfcup(\Omega\wedge(X_2\righthalfcup\Omega))-(X_1\righthalfcup X_2\righthalfcup\Omega)\Omega.$
 As the $2-$form $\Omega$ is characterising form and $X_2\righthalfcup\Omega$ is a linear combination of
 $\theta$ and $\phi$ , so that $\Omega\wedge(X_2\righthalfcup\Omega)=0$ implies that
 $(X_1\righthalfcup\Omega)\wedge(X_2\righthalfcup\Omega)=-(X_1\righthalfcup X_2\righthalfcup\Omega)\Omega.$
 Thus $\omega^1\wedge \omega^2=\frac{\Omega}{\Omega(X_2,X_1)}.$
 Clearly the process is inductive, and then the proof is complete.

To prove the second part, let \begin{eqnarray*}
 (\Omega(X_2,X_1))^2 d\bigg{(}\frac{\Omega}{\Omega(X_2,X_1)}\bigg{)}&=& (\Omega(X_2,X_1))^2(-\frac{d\Omega(X_2,X_1)
 \wedge \Omega}{(\Omega(X_2,X_1))^2}+\frac{d\Omega}{\Omega(X_2,X_1)})\\
 &=& -d\Omega(X_2,X_1)\wedge \Omega+\Omega(X_2,X_1)d\Omega\\
 &=& d(X_2\righthalfcup X_1\righthalfcup \Omega)\wedge\Omega+\Omega(X_2,X_1)d \Omega
  \end{eqnarray*}
Using $\pounds_{X_2} \alpha=
X_2\righthalfcup d\alpha+d(X_2\righthalfcup \alpha),$ we have
  \begin{eqnarray*}
 (\Omega(X_2,X_1))^2 d\bigg{(}\frac{\Omega}{\Omega(X_2,X_1)}\bigg{)} &=& [\pounds_{X_2}(X_1\righthalfcup \Omega)-X_2\righthalfcup d(X_1\righthalfcup\Omega)]\wedge\Omega
  +\Omega(X_2,X_1)d \Omega\\
   &=&-(X_2\righthalfcup (\pounds_{X_1}\Omega)\wedge\Omega+(X_2\righthalfcup X_1\righthalfcup d\Omega)\wedge\Omega
   +\Omega(X_2,X_1)d \Omega\\
   &=& -X_1\righthalfcup (X_2\righthalfcup d\Omega)\wedge\Omega+(X_2\righthalfcup d \Omega)
   \wedge(X_1\righthalfcup \Omega) +\Omega(X_2,X_1)d \Omega\\
\end{eqnarray*}
Now use Proposition \ref{prop:d(a^b)=k^a^b}, 
\begin{eqnarray*}
 (\Omega(X_2,X_1))^2 d\bigg{(}\frac{\Omega}{\Omega(X_2,X_1)}\bigg{)}
 &=& -X_1\righthalfcup \big{[}X_2\righthalfcup (K\wedge(X_1\righthalfcup\Omega)\wedge(X_2\righthalfcup
   \Omega))\big{]}\wedge\Omega+\Omega(X_2,X_1)d \Omega\\&+&\big{[}X_2\righthalfcup (K\wedge(X_1\righthalfcup\Omega)\wedge(X_2\righthalfcup
  \Omega))\big{]}\wedge(X_1\righthalfcup \Omega)\\
   &=& -\Omega(X_2,X_1)d\Omega+\Omega(X_2,X_1)d\Omega\\
   &=&0.
    \end{eqnarray*}
Hence, $d(\omega^1\wedge\omega^2)=0=d\bigg{(}\frac{\Omega}{\Omega(X_2,X_1)}\bigg{)}.$ The proof then follows by induction.
    \end{proof}
    \end{theorem}

The paper by Sherring and Prince \cite{sherr} extends Lie's approach to integrating a Frobenius integrable distribution by using solvable structures. The important results of \cite{sherr} are reproduced below. We illustrate these results with the straightforward case of a single second order ordinary differential equation.

 On some open subset of the first order jet bundle $J^1(\R,\R)$ with coordinates $x,y,p=\frac{dy}{dx},$ the second order differential equation (vector) field of the differential equation $\frac{d^2y}{dx^2}=f(x,y,p)$ is
\begin{equation}\label{soode}
\Gamma:=\frac{\partial}{\partial x}+p \frac{\partial}{\partial y}+f\frac{\partial}{\partial p}.
\end{equation}
This has dual distribution spanned by the contact form $\theta=dy-pdx$ and the force form $\phi=dp-dx.$ The 2-form $\Omega=\theta\wedge\phi$
on $U$ is then a characterising form for \eqref{soode}.

\begin{proposition}\cite{sherr}\label{Prop3}
  If $X$ is a vector field whose span with $\Gamma$ is two dimensional on $U$ then the $1-$form $X\righthalfcup \Omega$ is Frobenius integrable if and only if $X$ and $\Gamma$ are involution.
\end{proposition}

\begin{corollary}\cite{sherr}
If $X$ is a non-trivial symmetry of \eqref{soode}  then $X\righthalfcup \Omega$ is Frobenius integrable.
\end{corollary}

\begin{corollary}\cite{sherr}\label{Cor3.3 Sherr}
If $X_1$ and $X_2$ are symmetries of \eqref{soode} which commute and
whose span with $\Gamma$ is three dimensional then the two forms
$$\frac{X_1\righthalfcup\Omega}{X_2\righthalfcup X_1 \righthalfcup \Omega} \hspace*{15mm}
 {and\rm} \hspace*{15mm}\frac{X_2\righthalfcup\Omega}{X_1\righthalfcup X_2 \righthalfcup \Omega}$$
 are closed and locally  provide a complete set of two functionally independent first integrals.
\end{corollary}

The two first integrals of this corollary, together with $x,$ provide a coordinate chart which straightens out the integral curves of $\Gamma.$

The next corollary to Proposition~\ref{Prop3} shows how a solvable structure of two symmetries produces one closed form $\omega^1$ and a second form $\omega^2$ which is closed modulo $\omega^1,$ that is $d\omega^2\equiv0~ \text{mod}~\omega^1$ or locally $\omega^1=d\gamma^1$ and $\omega^2=d\gamma^2+\gamma^0d\gamma^1.$ This result in \cite{sherr} is proved by using the Frobenius theory of foliations but we will prove it in a very simple way by using exterior calculus. We need Theorem~\ref{thm2} and the following lemma.

\begin{lemma}\label{lemma:g mod f}
   Let $\alpha, \beta, d f\in \bigwedge^1(M)$ be non-zero one-forms with $\beta\nin\Sp\{df\}.$ Then $d\alpha=df\wedge\beta$ if and only if $\beta\equiv dg~\text{mod}~df~\text{for some}~g\in\bigwedge^0(M).$ Equivalently, $d\beta=\lambda\wedge df \iff \beta\equiv dg~\text{mod}~df.$

   \begin{proof}
Suppose that $d\alpha=df\wedge\beta.$ Then it is easy to see that $D^\perp:=\Sp\{df,\beta\}$ is Frobenius integrable and so locally there exists $dh$ such that $D^\perp:=\Sp\{df,dh\}$. Hence, $\beta=m~df+k~dh, ~m,~k\in \bigwedge^{0}$ and $d\alpha=df\wedge \beta=df\wedge(m~df+k~dh)=df\wedge k~dh.$ Now $0=d^2\alpha=d(df\wedge kdh)=-df\wedge dk \wedge dh,$ so that $dk\in \Sp\{df,dh\}$. This implies $k$ is a composite of $f$ and $h$ alone.
Let $g$ be the partial integral of $k$ with respect to $h$, then $k~dh=\pd{g}{h} dh$ and $d\alpha=df\wedge k~dh=df\wedge(\pd{g}{h} dh)=df\wedge(\pd{g}{h} dh+\pd{g}{f}df)=df\wedge dg.$
Hence, $df\wedge \beta=df\wedge dg$ which implies that $df\wedge (\beta-dg)=0~\text{and}~\beta\equiv dg~\text{mod}~df~\text{or}~\beta=dg+l~df.$ The converse is trivial.
   \end{proof}
\end{lemma}

 \begin{corollary}\label{cor3.4 Sherr}
If $X_1$ is a symmetry of \eqref{soode}, $X_2$ is a symmetry of $\Sp\{X_1,\Gamma\}$ and $X_1,X_2$ and $\Gamma$ are linearly independent everywhere then $$\omega^1:=\frac{X_1\righthalfcup\Omega}{X_2\righthalfcup X_1 \righthalfcup \Omega}$$ is closed, while
$$\omega^2:=\frac{X_2\righthalfcup\Omega}{X_1\righthalfcup X_2 \righthalfcup \Omega}$$ is closed modulo $\omega^1.$ These two forms locally provide a complete set of two functionally independent first integrals, $\omega^1=d\gamma^1$ and $\omega^2=d\gamma^2-X_2(\gamma^2)d\gamma^1.$ By putting $\hat{X_2}=X_2-X_2(\gamma^2)X_1$ we then have the two commuting symmetries $X_1$ and $\hat{X_2}$ which provide two first integrals $\gamma^1$ and $\gamma^2$ via Corollary \ref{Cor3.3 Sherr}.
\begin{proof}

The kernel of $X_1\righthalfcup\Omega$ is spanned by $X_1$ and $\Gamma$
and $X_1\righthalfcup\Omega$ is Frobenius integrable, since $
Sp(\{X_1,\Gamma\})$ is closed under lie bracket. So, $\omega^1$ is
closed by Theorem \ref{thm:sher2.1}. To see that $\omega^2$ is closed modulo
$\omega^1$, we consider
$\omega^1\wedge\omega^2=\frac{\Omega}{\Omega(X_2,X_1)}$ and by using the
fact $d(\omega^1\wedge\omega^2)=0$ (see proof in Theorem \ref{thm2}), this
implies that $d(\omega^1\wedge
\omega^2)=0=d(d\gamma^1\wedge\omega^2)=d\omega^2\wedge d\gamma^1.$

If $d\omega^2\wedge d\gamma^1=0$, from Theorem $1$ we may write  $d\omega^2=\mu\wedge d\gamma^1$ for some $\mu\in \Lambda^0$, so we have from Lemma \ref{lemma:g mod f} that $\omega^2=d\gamma^2+ld\gamma^1.$
To find the value of $l$, contract this equation with $X_1$ and $X_2$, we get $X_1(\gamma^2)=1$ and $l=-X_2(\gamma^2).$ Hence, $\omega^2=d\gamma^2-X_2(\gamma^2)d\gamma^1.$
\end{proof}
\end{corollary}
Corollary \ref{cor3.4 Sherr} can be generalised to the following theorem.
 \begin{theorem}\cite{sherr}\label{sherr result}
 Let $\Omega$ be a $k-$form on a manifold $M,$ and let $\Sp(\{X_1,\ldots, X_k\})$ be a $k-$dimensional distribution on an open $U\subseteq M$ satisfying $X_i\righthalfcup \Omega\neq 0$ everywhere on $U.$ Further suppose that $\Sp(\{X_{j+1},\ldots, X_k\}\bigcup\text{ker}~\Omega)$ is integrable for some $j<k$ and that $X_i$ is a symmetry of $\Sp(\{X_{i+1},\ldots, X_{k}\}\bigcup \text{ker}~\Omega)~\text{for}~ i=1,\ldots, j.$

 Put $\sigma^i:=X_1\righthalfcup \ldots \righthalfcup \bar{X}_i\righthalfcup\ldots\righthalfcup X_k\righthalfcup\Omega,$
where $\bar{X}_i$ indicates that this argument is missing and $\omega^i:=\frac{\sigma^i}{X_i\righthalfcup \sigma^i}~ \text{for}~i=1,\ldots,k$
so that $\{\omega^1,\ldots, \omega^k\}$ is dual to $\{X_1,\ldots, X_k\}.$ Then $d \omega^1=0;~d \omega^2=0~\mod~\omega^1;~d \omega^3=0~\mod~\omega^1,\omega^2;~~~\ldots~; d \omega^j=0~\mod~\omega^1,\ldots,\omega^{j-1},$ so that locally
 \begin{align*}
   \omega^1&=d \gamma^1,\\
   \omega^2&=d \gamma^2-X_1(\gamma^2)d\gamma^1,\\
   \omega^3&=d \gamma^3-X_2(\gamma^3)d\gamma^2-(X_1(\gamma^3)
   -X_2(\gamma^3)X_1(\gamma^2))d\gamma^1,\\ \vdots\\
  \omega^{j}&=d \gamma^{j}~\text{mod}~d\gamma^1,\ldots, d\gamma^{j-1},
 \end{align*}
 for some $\gamma^1,\ldots,\gamma^j\in \bigwedge^0 T^* U.$ Also the system $\{\omega^{j+1},\ldots,\omega^k\}$ is integrable modulo $d\gamma^1,\ldots, d\gamma^j$ and locally $\Omega=\gamma^0 d\gamma^1\wedge d\gamma^2 \wedge \ldots \wedge d\gamma^j\wedge \omega^{j+1}\wedge\ldots \wedge \omega^k ~\text{for some}~\gamma^0\in \bigwedge^0(T^* U).$ Each $\gamma^i$ is uniquely defined up to the addition of arbitrary function of $\gamma^1,\ldots,\gamma^{i-1}.$
 \end{theorem}

Sherring \cite{sherrexact} gives an algorithm for the integration of these forms which we reproduce in Theorem~\ref{IEM} below. We repeatedly use this process in our treatment of PDEs in subsequent sections. The integration of an exact $1-$form is described in the following lemma.
\begin{lemma}
Given an exact 1-form $\omega:=\omega_i~dx^i$ on $\R^n$, if
$\sigma^i:=\int(\omega_i-\sum_{j=1}^{i-1} \sigma^j_{x_i})~dx^{i}$ and
$\gamma:=\sum_{i=1}^{n}\sigma^{i}$ then $\omega=d \gamma.$
\begin{proof}
 Let $\omega=d\gamma$ for some $\gamma$ as $\omega$ is exact and so $\omega_i~dx^i=\frac{\partial \gamma}{\partial
 x^i}~dx^i.$
 For $i=1$, $\gamma=\int{\omega_1~d
 x^1}+\phi^1(x^2,\ldots ,x^n)=\sigma^1+\phi^1(x^2,\ldots ,x^n)$ and then,
 $\gamma_{x^2}=\sigma^1_{x^2}+\phi^1_{x^2}=\omega_2$
 this implies that, $\phi^1=\int{(\omega_2-\sigma^{1}_{x^2})~d x^2}+\phi^2(x^3,\ldots ,x^n)$
 and $\gamma=\sigma^1+\sigma^{2}+\phi^2(x^3,\ldots, x^n).$
 Now assume this holds for $i=k$,   $\gamma=\sum_{i=1}^{k}
 \sigma^{i}+\phi^{k}(x^{k+1},\ldots, x^n).$  Then the result holds for $i=k+1$,   $\gamma=\sum_{i=1}^{k+1}
 \sigma^{i}+\phi^{k}(x^{k+2},\ldots, x^n).$ The claim then follows by induction and the $n^{th}$  claim proves the theorem as $\phi^{n}$ is just the constant of
 integration.
\end{proof}
\end{lemma}

\begin{theorem}\label{IEM}
Let $\omega\in  \bigwedge(\R^n)$ which is exact modulo
$d\gamma^{1},\ldots, d\gamma^{j}.$ Without loss of generality choose coordinates
$x^{l}:=\gamma^{l}, \quad for \quad l=1,\ldots, j\ \ \text{and}~~~
y^k:=x^k, k=j+1,\ldots, n.$ Then $\omega:=\omega_i ~dx^i$, put
$\sigma^i:=\int(\omega_i-\sum_{k=j+1}^{i-1} \sigma^{k}_{x_i})dx^{i}$
for $i=j+1,\ldots, n$ and $\gamma:=\sum_{i=j+1}^{n}\sigma^{i}.$ Then
$\omega\equiv d \gamma  \mod d\gamma^1\ldots d\gamma^j.$
\begin{proof}
Given that $\omega$ is a non zero $1-$form which is exact modulo $d\gamma^{1},\ldots, d\gamma^{j}.$ Clearly, by
choosing coordinates $x^{l}=\gamma^{l}, \ l=1,\ldots, j$;  we then have $\omega$ is exact. Now 
$\omega_i=\frac{\partial \gamma}{\partial x^i}.$
 for $i=j+1$, \[\gamma=\int{\omega_{j+1} d
 x^{j+1}}+\phi^{j+1}(x^{j+2},\ldots, x^{n})=\sigma^{j+1}+\phi^{j+1}(x^{j+2},\ldots, x^n)\] and then,
 $\gamma_{x^{j+2}}=\sigma^{j+1}_{x^{j+2}}+\phi^{j+1}_{x^{j+2}}=\omega_{j+2}$
 this implies that, \[\phi^{j+1}=\int{(\omega_{j+2}-\sigma^{j+1}_{x^{j+2}})d
 x^{j+2}}+\phi^{j+2}(x^{j+3},\ldots, x^n)\]
 and
 \[\gamma=\sigma^{j+1}+\sigma^{j+2}+\phi^{j+2}(x^{j+3},\ldots, x^n).\]
 Clearly, for $i=j+3$,   $\gamma=\sum_{k=j+1}^{j+3}
 \sigma^{k}+\phi^{j+4}(x^{j+4},\ldots, x^n)$ and
 if we continue this process,
 $\gamma=\sum_{i=j+1}^{n} \sigma^{i}$ as $\phi^{n}$ is just the constant of
 integration.
\end{proof}
\end{theorem}

For ordinary differential equations the vector field description is one-dimensional, and hence integrable by dimension. For PDEs the situation is different. And so having explained the integration of a Frobenius integrable distribution via solvable structures we now turn to the issue of locating integrable sub-distributions of non-integrable distributions.

We begin with a non-integrable distribution $D$ of constant dimension $p$ on an $n$-dimensional manifold $M$. Let $D^\perp$ be the $n-p$ dimensional co-distribution and let $D^*$ be a complementary distribution with $\X(M)=D\bigoplus D^*.$ We aim to find the largest integrable sub-distributions of $D$, of dimension $p-M$ say, which satisfy the appropriate independence condition.

This is done in stages by using the following algorithm. We will try to reduce the dimension of $D^\perp$ and correspondingly add more $1-$forms to the distribution $D$.
Then the integrability condition \eqref{eqn:FI} for the augmented distribution is divided into two parts, one is  algebraic and the other is differential
 \begin{align}&d\sigma^a\wedge\Omega_\sigma\wedge\Omega_\rho=0, \ \ \ a=1,\ldots, n-p\label{eqn:FI1}\\
&d\rho^\alpha \wedge\Omega_\sigma\wedge\Omega_\rho=0, \ \ \ \alpha=1,\ldots, M\le p.\label{eqn:FI2}\end{align}
Here $\Omega_{\sigma}:=\sigma^1\wedge\ldots\wedge\sigma^{n-p}$ is a characterising form for $D$ and $\rho^\alpha$  are $1-$forms in $D^{*\perp}$ with $\Omega_\rho:=\rho^1\wedge\dots\wedge\rho^M$ and $\Omega_\sigma\wedge \Omega_\rho\neq 0.$

Firstly, we will add a single $1-$form $\rho\in D^{*\perp} $ to $D^\perp$ and generate the algebraic conditions \eqref{eqn:FI1} for $\Sp\{\sigma^1,\ldots\sigma^p,\rho\}$ to be Frobenius integrable. If the algebraic conditions can be satisfied then there may be integrable Vessiot sub-distributions of dimension $p-1$. If the algebraic conditions fails then we will add two $1-$forms from $D^{*\perp}$. The process continues until we have an enlarged co-distribution which may be integrable, if $M=p-1$ then the result is trivial. The next step is to solve the differential conditions \eqref{eqn:FI2}. This is the difficult part and in general the best we can do is apply an EDS process which will identify the number and  functional dependency of the solutions rather than explicitly construct solutions. We explain this process in the last section. In the event that we can construct a closed form solution we can, in the presence of a solvable structure, integrate the Frobenius integrable sub-distribution by using Theorem \ref{sherr result}.

In the following section, we will consider the particular type of second order PDEs of one dependent variable and two independent variables and explain the algorithm in detail.

\section{Second order PDEs}
\noindent In this section, we examine two types of second order PDEs of one dependent variable and two independent variables. The first is of the form
$u_{yy}= F(x,y,u,u_x,u_y,u_{xx},u_{xy})$ and the second PDEs of evolution type.
\subsection{PDEs of the type $u_{yy}= F(x,y,u,u_x,u_y,u_{xx},u_{xy})$}
\noindent Consider a second order PDE with one dependent real variable $u$ and two independent real variables $x,y$ given by
 \begin{equation}\label{SOPDE}
 u_{yy}= F(x,y,u,u_x,u_y,u_{xx},u_{xy}),
 \end{equation}
 for some smooth function $F.$
The embedded submanifold \cite{Stormark} \[S:=\{(x,y,u,u_x,u_y,u_{xx},u_{xy})\in J^2(\mathbb{R}^2,\mathbb{R})\ | \ \ u_{yy}-F(x,y,u,u_x,u_y,u_{xx},u_{xy})=0\},\] is a subset of $J^2(\mathbb{R}^2,\mathbb{R}).$
 A local solution of the PDE is a  $7-$dimensional locus of $J^2(\mathbb{R}^2,\mathbb{R})$ described by the map $i: S \rightarrow J^2(\mathbb{R}^2,\mathbb{R}),$ i.e. \[i:(x,y,u,u_x,u_y,u_{xx},u_{xy})\mapsto (x,y,u,u_x,u_y,u_{xx},u_{xy},F).\]

We can study the solutions of \eqref{SOPDE}  by studying the integral submanifolds $N$                                                                                                                                                                                                                                                                                                                                                                                                                                                                                                                                                                                                                                                                                                                                                                                                                 of the restricted contact system
\begin{align}
D^{\perp}_V:= \Sp \big{\{}\theta^1&:=du-u_xdx-u_ydy,\nonumber\\
 \theta^2&:=du_x-u_{xx}dx-u_{xy}dy,\label{contact dist} \\
 \theta^3&:=du_y-u_{xy}d x-Fdy\big{\}}\nonumber,
 \end{align}
which project down to $X\subset \mathbb{R}^2.$
If $N$ satisfies the transverse condition $dx\wedge dy| _N \neq 0$ and has a tangent space that annihilates the restricted contact structure, then $i(N)\subset J^2(\mathbb{R}^2,\mathbb{R})$ is the $2-graph$ of a solution of \eqref{SOPDE}. So, we need to study the vector field distribution dual of the pulled-back contact system, that is,
\begin{align}
D_V:= \Sp\bigg{\{V_1}&:=\vp{x}+u_{x}\vp{u}+u_{{x}{x}}\vp{u_{x}}+u_{xy}\vp{u_{y}},\nonumber\\
V_2&:=\vp{y}+u_{y}\vp{u}+u_{xy}\vp{u_{x}}+F\vp{u_{y}},\label{Vessiot dist}\\
V_3&:=\vp{u_{{x}{x}}},V_4:=\vp{u_{xy}}\bigg{\}}.\nonumber
\end{align}
This dual distribution is called the  Vessiot distribution $D_V$ of \eqref{SOPDE}.

\subsubsection{Algorithm for finding the largest integrable sub-distributions of $D_V$}\quad

\noindent In general the Frobenius integrability condition fails on the Vessiot distribution $D_V$. So, we will try to reduce the dimension of $D_V$ and correspondingly add more $1-$forms from $D^{\star\perp}$ to the co-distribution $D^{\perp}_V$ as described at the end of the Section $3$.

For this class of PDEs we need to add exactly two $1-$forms to $D_V^\perp$ because in adding a single $1-$form $\phi$ the second algebraic condition $d\theta^2\wedge\Omega_\theta\wedge\phi=0$ implies that $\phi=0.$ Furthermore, we cannot add more than two $1-$forms because the reduced Vessiot distribution cannot have dimension less than 2.

\noindent \textbf{Step 1.}
Add two 1-forms $\phi^1, \phi^2\in D^{\star\perp} $ in $D_V^{\perp}.$

\noindent We are interested in the largest integrable sub-distributions which satisfy the independence condition ($d x\wedge dy\neq0$). So, for this purpose the coefficient of $\psi^3$ and $\psi^4$ in $\phi^1,~\phi^2$ should be nonzero. Without a loss of generality we assume that
\begin{equation}\label{eqn:phi1,2}
\phi^1:=\psi^3-a_{1} \psi^1-a_{2}\psi^2,\ \ \phi^2:=\psi^4-b_{1}\psi^1-b_{2}\psi^2.
\end{equation}
 Now, solve the algebraic conditions \[d\theta^a\wedge\Omega_{\theta}\wedge\Omega_\phi=0,\ \ a=1\ldots 3.\]
 Specifically, \begin{align*}
  a_2=b_1,\quad \quad b_2=u_xF_u+u_{xx}F_{u_x}+a_1F_{u_{xx}}+b_1F_{u_{xy}}+u_{xy}F_{u_y}+F_x.
\end{align*}

\noindent\textbf{Step 2.} Differential conditions:\label{step:2diff cond}

\noindent Solve the differential conditions \eqref{eqn:FI2} i.e.
\begin{align*}
  d\phi^1 \wedge\Omega_\theta\wedge\Omega_\phi=0,\\
  d\phi^2 \wedge\Omega_\theta\wedge\Omega_\phi=0.
\end{align*}

\noindent\textbf{Step 3.}\label{step3:integration} Integrate the reduced Vessiot distribution $D_{V_{red}}$ by using  Theorem \ref{sherr result}.
To find the symmetries of reduced Vessiot distribution use DIMSYM, for example.

\begin{remark}
In applying the algorithm, there will exist situations when it may be difficult to solve the differential conditions (step $2$). For such situations we will look on the alternative way to solve them. In Section $5,$ we will discuss the solution of these differential conditions in detail.
\end{remark}
We illustrate this algorithm with an example:
\begin{example}
Consider the Laplace equation
$$u_{xx}+u_{yy}=0$$
The Vessiot distribution and co-distribution is given by equations \eqref{contact dist} and \eqref{Vessiot dist}  with $F=-u_{xx}.$\\
Now apply the algorithm step by step. By adding the two $1-$forms \eqref{eqn:phi1,2} in the $D^{\perp}_V$ and after solving the algebraic conditions \eqref{eqn:FI1} we obtain $a_{2}  =b_1,~~a_{1}  = - b_{2}.$

The corresponding differential constraints are
\begin{align*}
 d\phi^1\wedge\Omega_\theta\wedge\Omega_\phi&=u_x\pd{b_1}{u}-u_y\pd{b_2}{u}-\pd{b_1}{x}-\pd{b_2}{y}-u_{xx}\pd{b_1}{u_x}-u_{xy}\pd{b_2}{u_{x}}-u_{xy}\pd{b_1}{u_y}+u_{xx}\pd{b_2}{u_y}\\
 &+b_2\pd{b_1}{u_{xx}}-b_1\pd{b_2}{u_{xx}}-b_1\pd{b_1}{u_{xy}}-b_2\pd{b_2}{u_{xy}},\\
 d\phi^2\wedge\Omega_\theta\wedge\Omega_\phi&=u_y\pd{b_1}{u}-u_x\pd{b_2}{u}-\pd{b_2}{x}+\pd{b_1}{y}+u_{xy}\pd{b_1}{u_x}-u_{xx}\pd{b_2}{u_{x}}-u_{xx}\pd{b_1}{u_y}-u_{xy}\pd{b_2}{u_y}\\
 &+b_1\pd{b_1}{u_{xx}}+b_2\pd{b_2}{u_{xx}}-b_2\pd{b_1}{u_{xy}}-b_1\pd{b_2}{u_{xy}}.
\end{align*}In order to show some solutions we will choose an arbitrary value for these coordinates.
For a particular choice of parameters we have the following solution.

$\phi^1:=\psi^3-\psi^2, \phi^2:=\psi^4-\psi^1$ satisfy the differential conditions and the transverse condition. The corresponding maximal dimension integrable sub-distribution  $D_{red}$ is generated by
\begin{align*}
\bar{V}_1&:=\vp{x}+u_x\vp{u}+u_{xx}\vp{u_x}+u_{xy}\vp{u_y}+\vp{u_{xy}},\\
\bar{V}_2&:=\vp{y}+u_y\vp{u}+u_{xy}\vp{u_x}-u_{xx}\vp{u_y}+\vp{u_{xx}}.
\end{align*}
For the integration of reduced Vessiot distribution use Theorem \ref{sherr result} and the following symmetries of $D_\text{red}$
\begin{align*}
X_1&:=\vp{u},\quad \quad X_2:=-x\vp{u}-\vp{u_x},\quad \quad X_3:=-y\vp{u}-\vp{u_y},\\
X_4&:=-xy\vp{u}-y\vp{u_x}-x\vp{u_y}-\vp{u_{xy}},\quad X_5:=-(y^2-x^2)\vp{u}-2y\vp{u_y}-2x\vp{u_x}-2\vp{u_{xx}},
\end{align*}
and after integration we obtain the following invariant functions
\begin{align*}
f^1&:=u-xu_x+xyu_{xy}-yu_y+\frac{1}{2}(u_{xx}x^2-u_{xx}y^2-x^2y)+\frac{1}{6}y^3,\\
f^2&:=u_{x}-yu_{xy}-xu_{xx}+xy,\quad \quad f^3:=x u_{xy}-y u_{xx}-u_y+\frac{1}{2}(y^2-x^2),\\
f^4&:=u_{xy}-x,\quad \quad f^5:=\frac{1}{2}(u_{xx}-y).
\end{align*}

\noindent The lifted solution on $S$ is a common level set of these functions:
$$\left\{ p\in S:\ f^\alpha(p)=c_\alpha\right\}.$$
This projects to
$u=c_1+c_2x- c_3y+c_4xy+\frac{1}{2}c_5(x^2-y^2)+\frac{1}{2}x^2y-\frac{1}{6}y^3.$
\end{example}

\begin{example}
Consider the nonlinear second order partial differential equation (generalised equation of steady transonic gas flow) \cite{Polyanin}
\[u_{yy}+A\frac{u_y}{y}+Bu_xu_{xx}=0,\quad y\neq 0,\]
where $A$ and $B$ are arbitrary constants. For the sake of simplicity we choose $A=1$ and $B=1.$ Now apply the algorithm step by step.
By adding the two $1-$forms \eqref{eqn:phi1,2} in the $D^{\perp}_V$ and the solution of the algebraic conditions \eqref{eqn:FI1} implies
$a_{2}  =b_1,~~b_{2}  = a_1u_x+u_{xx}^2+\frac{u_{xy}}{y}.$

For a particular choice of parameters we find that
$$\phi^1:=\psi^3,\quad \phi^2:=\psi^4+(u^2_{xx}+\frac{u_{xy}}{y})dy$$  satisfy the differential conditions (we suppress domain considerations for clarity). The corresponding maximal dimension integrable sub-distribution  $D_{red}$ is generated by
\begin{align*}
\bar{V}_1&:=\vp{x}+u_x\vp{u}+u_{xx}\vp{u_x}+u_{xy}\vp{u_y},\\
\bar{V}_2&:=\vp{y}+u_y\vp{u}+u_{xy}\vp{u_x}+(u_xu_{xx}+\frac{u_y}{y})\vp{u_y}-(u^2_{xx}+\frac{u_{xy}}{y})\vp{u_{xy}}.
\end{align*}
We now proceed to apply Theorem \ref{sherr result} to integrate the reduced Vessiot distribution. The symmetries of $D_\text{red}$ are
\begin{align*}
X_1&:=\vp{x},\quad \quad X_2:=2x\vp{x}+y\vp{y}+4u\vp{u}+2u_x\vp{u_x}+3u_y\vp{u_y}+u_{xy}\vp{u_{xy}},\quad \quad X_3:=\vp{u},\\
X_4&:=\log y\vp{u}+\frac{1}{y}\vp{u_{y}},\quad X_5:=2u\vp{u}-y\vp{y}+2u_x\vp{u_x}+3u_y\vp{u_y}+3u_{xy}\vp{u_{xy}}+2u_{xx}\vp{u_{xx}}.
\end{align*}
After integration, we obtain the following conserved quantities:
\begin{align*}
f^1&:=\frac{1}{4u_{xx}}(4xu_{xx}-4u_x - y^2u_{xx}^2+(y^2u_{xx}^2+ 2yu_{xy})(\log u_{xx}+2\log y)),\\
f^2&:=\frac{1}{2}(-\log y+\log u_{xx}-\log (2u_{xy}+yu_{xx}^2),\\
f^3&:=\frac{1}{64u_{xx}}(16yu_yu_{xx}- 16yu_xu_{xy}- 3y^4u_{xx}^4- 8y^3u_{xx}^2u_{xy}+4\log y ),\\
f^4&:=\frac{1}{64u_{xx}}(32u_x^2+7y^4u_{xx}^4+16y^3u_{xx}^2u_{xy}-(\log y)^2( 8y^4u_{xx}^4+32y^3u_{xx}^2u_{xy}+32y^2u_{xy}^2-16)\\
&-\log y( 64yu_xu_{xy}+12y^4u_{xx}^4+32y^3u_{xx}^2u_{xy}-16yu_yu_{xx})),\\
f^5&:=\frac{1}{2}\log u_{xx}.
\end{align*}
As before we can express the solution in the original coordinates.
\end{example}

\subsection{Evolution equations}
\noindent In this section we present an additive separable solution of the evolution equations. Consider the following second order evolution equation
\begin{equation}
u_{x}= F(x,y,u,u_y,u_{yy}),
\end{equation} with an extra condition $u_{xy}=0.$ This extra condition means the solution is additively separable. In \cite{Vassiliou98} the authors say that a function has an additively separable form $u(x,y)=v(x)+w(y)$ if and only it satisfies the condition $u_{xy}=0.$

The restricted contact distribution on $\R^5$ is
\begin{align*}
D^{\perp}_V:= \Sp \big{\{}\theta^1&:=du-Fdx-u_ydy,\\
 \theta^2&:=du_y-u_{yy}dy\big{\}}.
 \end{align*}
The corresponding dual distribution is
\begin{align*}
D_V:= \Sp\bigg{\{V_1}&:=\vp{x}+F\vp{u},~~
V_2:=\vp{y}+u_{y}\vp{u}+u_{yy}\vp{u_{y}},~~
V_3:=\vp{u_{yy}}\bigg{\}}.
\end{align*}
In this case $D^{*\perp}:=\Sp\left\{\psi^1=dx, \psi^2=dy, \psi^3=du_{yy} \right\}.$

Now follow the steps of the above algorithm:\\
If we add a 1-form $\phi\in D^{\star\perp} $ to $D_V^{\perp}$ and generate the algebraic conditions \eqref{eqn:FI1} for $\Sp\{\theta^1,\theta^2,\phi\}$ to be Frobenius integrable we find a  non zero solution. So, there is an integrable Vessiot sub-distribution of dimension 2.
The solution is
\begin{align}\label{eqn:evol alg sol}
   a_1&=0,\nonumber\\
  a_2&=-\frac{F_y+u_yF_{u}+u_{yy}F_{u_y}}{F_{u_{yy}}}, \quad F_{u_{yy}}\neq0.
\end{align}
For the evolution equations we have found the solution of algebraic conditions by adding one $\phi:=\psi^3-a_1\psi^1-a_2\psi^2$ in the $D_V^{\perp}.$ So the process is very simple in this case we need to go to the step $3$ directly which is a differential condition. At this stage we need an expression for $F.$ For this purpose we consider the potential Burgers equation.\newpage
\begin{example}
Consider the potential Burgers equation
  \[ u_x=u_{yy}+u_y^2,\]
 with $u_{xy}=0.$ From \eqref{eqn:evol alg sol}, we have  $a_1=0,~a_2=-2u_y.$ This implies
 \[\phi:=\psi^3+2u_y\psi^2\] which satisfy the differential condition and transverse condition.\newline
The reduced Vessiot distribution is spanned by
\[D_\text{red}:=\Sp\{\bar{V}_1:=V_1,~~~~~\bar{V}_2:=V_2-2u_yu_{yy}V_3\}.\]
The final step is the integration of reduced Vessiot distribution using Theorem \ref{sherr result}.

The symmetries of $D_\text{red}$ are
\begin{align*}
X_1&:=\vp{y},\quad \quad X_2:=\vp{u},\quad \quad X_3:=2x\vp{x}+y\vp{y}-u_y\vp{u_y}-2u_{yy}\vp{u_{yy}},
\end{align*}
and after integration, we obtain
\begin{align*}
f^1&:=y-\frac{1}{2}(u_y^2+u_{yy})^{-\frac{1}{2}}\text{arctanh}((u_y^2+u_{yy})^{-\frac{1}{2}}u_y),\\
f^2&:=u-(u_y^2+u_{yy})x+\frac{1}{2}\ln(u_{yy})-\frac{1}{2}\ln(u_y^2+u_{yy}),\\
f^3&:=\frac{1}{2}\ln(u_y^2+u_{yy}).
\end{align*}
The solution is $u=c_2+xe^{2c_3}-\frac{1}{2}\ln{(1-\tanh^2{(ye^{c_3}-c_1e^{c_3})})}.$

\end{example}
\begin{remark}
It is worth pointing out the utility of our geometric approach in the investigation of some natural questions concerning conserved quantities and symmetries. Having solved the problem of finding all the maximal integrable sub-distributions of the Vessiot distribution we can now ask, for example, whether there exists a solution of the PDE admitting a conserved quantity of energy type, or whether there exists a solution simultaneously admitting certain specified conserved quantities. We can also ask whether there exist solutions admitting a particular symmetry algebra, solutions which have certain topological properties and so on.

As an example, consider the Gaussian curvature of the graph of a solution of the PDE $u_{yy}=F(x,y,u,u_x,u_y,u_{xx},u_{xy})$:
$$K=\frac{u_{xx}F-u_{xy}^2}{(1+u_x^2+u_y^2)^2}$$
and suppose that this solution is the projection of a common level set of five functions $f^a$ of the type we have discussed earlier.
Then $K$ is clearly constant on the graph of the solution if and only if $dK\wedge df^1\wedge \ldots\wedge df^5=0,$ expressing that $K$ is a composite of $f^1,\ldots,f^5.$
But less obvious and much more useful is that $K$ is constant on the graph if and only if $dK\wedge \Omega_{\theta}\wedge \Omega_{\phi}=0.$ Moreover, if $K$ is not constant on the graph then the one dimensional integral manifolds of $dK\wedge \Omega_{\theta}\wedge \Omega_{\phi}$ are the curves on the graph along which $K$ is constant and their local existence is guaranteed.\newline
Further calculus with $K$ on the solution submanifold of $J^2(\R^2,\R)$ is effected by using the tangent vector fields $\bar{V}_1,~\bar{V}_2.$ Since $\bar{V}_1,~\bar{V}_2$ project to the tangent directions to the graph of the solution, this calculus also projects. The point being that calculus on $J^2(\R^2,\R)$ with functions like $K$ is easier than on the base graph space.
\end{remark}

\section{Solving the differential conditions}
\noindent In this section we discuss the solution of differential conditions \eqref{eqn:FI2}. We present two methods. The first is an Exterior Differential System (EDS) method and the second builds a solution invariant under a solvable structure.
\subsection{EDS method}\quad

\noindent From Proposition \ref{prop:d(a^b)=k^a^b}, we know that solving the Frobenius conditions
 \begin{align*}&d\theta^a\wedge\Omega_\theta\wedge\Omega_\phi=0, \ \ \ a=1,\ldots, 3\\
&d\phi^\alpha \wedge\Omega_\theta\wedge\Omega_\phi=0, \ \ \ \alpha=1,\ldots, M\leq 2\end{align*}
for $\Omega_\phi:=\phi^1\wedge\dots\wedge\phi^M$ is equivalent to solving
\begin{equation}d(\Omega_\theta \wedge \Omega_\phi)=\lambda\wedge\Omega_\theta\wedge \Omega_\phi \label{eqn:dif}\end{equation}
for $\Omega_\phi$ and some 1-form $\lambda.$ For later reference we note that, without loss, $\lambda$ has no $\theta^a$ or $\phi^\alpha$ components, and that  $d\lambda\wedge\Omega_\theta\wedge\Omega_\phi=0.$ This last condition means that $d\lambda$ only has components with $\theta^a$ or $\phi^\alpha$ as a factor.

Our approach here will be to create a differential ideal of 5-forms and then, through a closure condition, find those $\Omega_\theta\wedge\Omega_\phi$ which satisfy \eqref{eqn:dif}.

\noindent The standard EDS reference are the books \cite{Ivey,bryant} and more details of EDS method by means of the inverse problem can be found in the work of \cite{Anderson1992}. We will give a brief summary of the method in this context.

The EDS process for finding the 5-forms satisfying \eqref{eqn:dif} involves three steps. For the case of second order partial differential equations, we start with a module of 5-forms, $\Sigma$, and find the submodule of $\Sigma^\text{final}$ that forms a differential ideal. We will find (or not) our Frobenius integrable  5-form in this ideal. The second step is to reformulate the Frobenius condition in terms of an equivalent linear Pfaffian system of one forms, and final step is to use the Cartan-K\"{a}hler theorem to determine the generality of the solutions of the problem.

The actual process for finding the  submodule $\Sigma^\text{final} \subset \Sigma$ that generates the differential ideal is the following terminating recursive procedure. Starting with the submodule $\Sigma^0:=\Sigma$ find the submodule $\Sigma^1\subseteq\Sigma^0$ such that for all $\Omega\in \Sigma^1,$ $d\Omega\in\langle\Sigma^0\rangle$, the ideal generated by $\Sigma^0.$

We check if $\Sigma^1=\Sigma^0$ and so is already a differential ideal. If not, we iterate the process, finding the submodule $\Sigma^2\subset\Sigma^1\subset\Sigma^0$ and so on until at some step, a differential ideal is found or the empty set is reached. If at any point during the process it is not  possible to find a submodule containing simple forms, the problem has no solution.

Having found $\Sigma^\text{final},$ the next step in the EDS process is to express the problem of finding the Frobenius integrable 5-form in $\Sigma^\text{final}$ as a Pfaffian system. Let the submodule $\Sigma^\text{final},$ be spanned by the set $\{\Omega^k:k=1,\dots,d\},$ and calculate \[d\Omega^k=\xi^k_h\wedge\Omega^h\] to find the $\xi^k_h$'s which are now fixed one forms.

Let $\Omega=r_k\Omega^k,$
now we are looking for those $\Omega$'s in $\Sigma^\text{final}$ such that $d\Omega=\lambda\wedge\Omega,$ this gives:
\begin{align}
  d(r_k\Omega^k)=\lambda\wedge r_k\Omega^k\nonumber,\\
  \Leftrightarrow dr_k\wedge\Omega^k+r_kd\Omega^k-r_k\lambda\wedge\Omega^k=0,\nonumber\\
  \Leftrightarrow (dr_k+r_h\xi^h_k-r_k\lambda)\wedge\Omega^k=0. \label{eqn:dr_k1}
\end{align}
So this leads us to the task of finding all the possible $d-$tuples of one forms $(\rho_k)=(\rho_1,...,\rho_d)$ such that
\[\rho_k\wedge\Omega^k=0.\] Suppose the solutions are
\[(\rho^A_k):=(\rho_1^A,...,\rho^A_d),~A=1,...,e\] with
\[\rho_k^A\wedge\Omega^k=0,\] that is, an $e-$dimensional module of $d-$tuples of $1-$forms. Then equation (\ref{eqn:dr_k1}) can be expressed as
\begin{align}
  (dr_k+r_h\xi^h_k-r_k\lambda)=-p_A\rho^A_k,\nonumber \\
 \Leftrightarrow dr_k+r_h\xi^h_k-r_k\lambda+p_A\rho^A_k=0,\label{eqn:dr_k}
\end{align}
where $p_A$ are arbitrary functions.

At this stage, firstly the problem becomes that of solving (\ref{eqn:dr_k}) for $r_k$ in terms of the unknown functions $p_A,$ and secondly, identifying the restrictions on the choice of these $p_A$'s.
The general method for finding the solution for this problem in EDS is to define an extended manifold $N:=S\bigotimes\R^d\bigotimes\R^e$ with coordinates $\{x,y,u,u_x,u_y,u_{xx},u_{xy},r_k,p_A\}, k\in \{1,...,d\}, A\in\{1,...,e\}$ and look for $7-$dimensional manifolds that are sections over $S$ and on which the one forms
\begin{equation}\sigma_k:=dr_k+r_h\xi^h_k-r_k\lambda+p_A\rho^A_k\label{eqn:sigma}\end{equation} are zero.
To find these manifolds, $\sigma_k$ are considered constraint forms for a distribution on $N$ whose integral submanifolds we want. We choose a basis of forms on $N,$ $\{\alpha_a:=dx,dy,du,du_{x},du_y,du_{xx},du_{xy},\sigma_k,\pi_A\}$ where $\alpha_a$ are pulled-back basis for $S,$ $\pi_A:=dp_A$ and $\sigma_k$ are defined above completes the basis. The condition that we want sections over $S$ is that $\alpha_1\wedge \alpha_2\neq 0.$
The Frobenius integrability condition for $D_\sigma=\Sp\{\sigma_k\}$ is $\d\sigma_k\equiv 0\mod\sigma_k.$
But,\begin{equation}d\sigma_k\equiv\pi^i_k\wedge\alpha_i+t^{12}_k\alpha_1\wedge\alpha_2~~\mod~~\Sp\{\sigma_k\}\label{eqn:dsigma}\end{equation} where the $\alpha_1\wedge\alpha_2$ part is the \emph{torsion} and $\pi^i_k$ are some linear combination of $dp_A.$ As $d\sigma_k$ expands with no $dp_A\wedge dp_B$ terms, the system is quasi-linear.

Because we want $\alpha_1\wedge\alpha_{2}\neq 0$ on the integral manifolds, we need to absorb all the \emph{torsion} terms into the $\bar{\pi}_A=\pi_A-l^j_A\alpha_j.$ If any of the \emph{torsion} terms cannot be absorbed, then asking for $d\sigma_k\equiv0~~\mod~~\Sp\{\sigma_k\}$ is incompatible with the independence condition and therefore there is no solution.

Once the \emph{torsion} terms have been removed, so that (\ref{eqn:dsigma}) becomes
\begin{equation}
  d\sigma_k\equiv\bar{\pi}^i_k\wedge\alpha_i,\label{eqn:dsigbar}
\end{equation}
the next step is to create the tableau $\prod$ from (\ref{eqn:dsigbar}), from which the Cartan characters can be calculated allowing us to apply the Cartan test for involution.

In order to satisfy the integrability condition we will force all the $\bar{\pi}^i_k$ terms to be zero. As a result we will obtain an over-determined system of linear algebraic equations for $r_k.$

We illustrate this idea with an example:
\begin{example}
\noindent Consider the second order Laplace equation $u_{yy}=-u_{xx}.$ In order to produce an initial ideal we solve the algebraic equation \eqref{eqn:FI1} giving
\begin{align*}
 a_2=b_1,~a_1=-b_2.
\end{align*}
This means that $\Omega_\phi$ is in the submodule, $\Sigma^0$, of $2-$forms  generated by
\begin{align*}
  \omega^1&:=-\psi^1\wedge\psi^2,\\
  \omega^2&:=\psi^1\wedge\psi^4+\psi^2\wedge\psi^3,\\
  \omega^3&:=\psi^1\wedge\psi^3-\psi^2\wedge\psi^4,\\
  \omega^4&:=\psi^3\wedge\psi^4.
\end{align*}
Following the EDS process above, the first step is to find a differential ideal. This is done in an iterative process as follows.

 Starting with $\Sigma^0=\Sp\{\Omega^k\},~~ k=1,...,4,$ where $\Omega^k=\omega^k\wedge\Omega_\theta,$ find the submodule $\Sigma^1$ where for each $\bar{\Omega}\in \Sigma^1,~d\bar{\Omega}$ is in $\langle\Sigma^0\rangle,$ the ideal generated by $\Sigma^0. $ In this example we have $\Sigma^1=\Sigma^0.$ It means we already have a differential ideal and the process terminates.

 So, now having  a differential ideal of dimension four, we can continue with the EDS process. The next step is to solve $d(\omega^k\wedge\Omega_\theta)=\xi^k_l\wedge\omega^l\wedge\Omega_\theta,$ for the one forms $\xi^k_l.$ For this particular example we find that all
$\xi^k_l=0.$\\
The next step is to find $4-$tuples of one forms $(\rho_1,...,\rho_4)$ such that:
\begin{equation}
  \rho_k\wedge\omega^k\wedge\Omega_\theta=0,\label{eqn:rho}
\end{equation}
 Let $  \rho_k=\hat{\rho}_{k}+\tilde{\rho}_{k},$
with $\hat{\rho}_{k}=\hat{a}_{ki}\psi^i,\tilde{\rho}_{k}=\tilde{a}_{ka}\theta^a.$ In this case, the $16~\hat{a}_{ki}$ are real constants satisfying
\begin{align*}
  \hat{a}_{23}&=-\hat{a}_{34}-\hat{a}_{41},\quad \quad  \hat{a}_{13}=-\hat{a}_{21}-\hat{a}_{32},\quad \quad  \hat{a}_{31}=\hat{a}_{22}-\hat{a}_{14},\quad \quad  \hat{a}_{42}=-\hat{a}_{33}+\hat{a}_{24},
\end{align*}
and the $12~\tilde{a}_{ka}$ are arbitrary. So we have a $24$ parameter family of solutions $\rho^A_k,~~A=1,...,24,$ with $\hat{\rho}_k^A=0~\text{for}~A=13,\ldots, 24 ~\text{and}~\tilde{\rho}_k^A=0~\text{for}~A=1,\ldots, 12. $
Following the procedure from the outline, we now extend $S$ to a new manifold with coordinates: $(x,y,u,u_x,u_y,u_{xx},u_{xy},r_k,p_A)$ and the problem can be written as that of finding the integrable distributions on $N$ with $\sigma_k=0$ where
\begin{equation}
\sigma_k:=dr_k-r_k\lambda+p_A\rho^A_k,\label{eqn:sigma_k}
\end{equation}
along with independence condition $dx\wedge dy\neq 0$ on these distributions. Continuing the EDS process, set $\pi_A=dp_A,$ Using this, a co-frame on $N$ is $(dx,dy,du,du_x,du_y,du_{xx},du_{xy},\sigma_k,\pi_A).$ So the next step is to calculate $d\sigma_k$ modulo the ideal $\langle\sigma_k\rangle$ as follows:

 \noindent Taking the exterior derivative of \eqref{eqn:sigma_k} and using this equation to remove resulting appearances of $dr_k$, we have
$$d\sigma_k\equiv (\pi_A-p_A\lambda)\wedge\rho^A_k-r_kd\lambda+p_Ad\rho^A_k,\quad \quad \quad  \mod\{\sigma_k\}$$
As we know that we have a $24$ parameter family of solutions, so we will consider the precise expressions for $\lambda$ and $d\lambda$ for further simplification. Recall that $\lambda$ has no $\theta^a$ components and that $d\lambda$ has only $\theta\wedge\psi$ and $\psi\wedge\psi$ components.
Let\begin{align*}
    \lambda:=\lambda_{i}\psi^i\ \quad \text{and}\quad d\lambda_i:=\lambda_{ij}\psi^j + \tilde{\lambda}_{ia}\theta^a.
  \end{align*}
\begin{align*}
  d\sigma_1&\equiv(\pi_1-p_1\lambda-r_1d\lambda_1)\wedge\psi^1+(\pi_2-p_2\lambda-r_2d\lambda_2)\wedge\psi^2
  +(\pi_4+\pi_7-(p_4+p_7)\lambda-r_1d\lambda_3)\wedge\psi^3\\&+(\pi_3-p_3\lambda-r_1d\lambda_4)\wedge\psi^4+\sum_{A=13}^{24}[(\pi_A-p_A\lambda)\wedge\rho^A_1+p_Ad\rho^A_1]\ \mod~~\{\sigma_k\}.
\end{align*}
In this case there exist $r_k$ with $dr_k\in\Sp\{\psi^i\}$ such that $$\lambda=(2-u_{xx}\lambda_3-u_{xy}\lambda_4)\psi^1+(u_{xx}\lambda_4-u_{xy}\lambda_3)\psi^2+\lambda_3\psi^3+\lambda_4\psi^4,~\text{where}~\lambda_3,~\lambda_4~ \text{are arbitrary,}$$ which can be seen by inspection from \eqref{eqn:sigma_k} with $p_A=0,~A=13,\ldots,24.$

One possible solution:
$$ r_1=u_{xx}^2+u_{xy}^2,~~~~r_2=-u_{xx},~~~~r_3=u_{xy},~~~~~r_4=1,$$ and
\begin{align*}
\Omega_\phi=-(u_{xx}^2+u_{xy}^2)\psi^1\wedge\psi^2-u_{xx}(\psi^1\wedge\psi^4+\psi^2\wedge\psi^3)+u_{xy}(\psi^1\wedge\psi^3-\psi^2\wedge\psi^4) +\psi^3\wedge\psi^4.
\end{align*}
The distribution $\Sp\{\theta^1,\theta^2,\theta^3,\phi^1,\phi^2\}$ is Frobenius integrable. We now proceed to apply the step $4$ of the algorithm.
The reduced Vessiot distribution is
\begin{align*}
\bar{V}_1&:=\vp{x}+u_x\vp{u}+u_{xx}\vp{u_x}+u_{xy}\vp{u_y}+u_{xx}\vp{u_{xx}}+u_{xy}\vp{u_{xy}},\\
\bar{V}_2&:=\vp{y}+u_y\vp{u}+u_{xy}\vp{u_x}-u_{xx}\vp{u_y}+u_{xy}\vp{u_{xx}}-u_{xx}\vp{u_{xy}}.
\end{align*}
After integrating the reduced Vessiot distribution, we have
\begin{align*}
f^1&:=e^{-x}(u_{xx}\cos y-u_{xy}\sin y),\quad \quad \quad \quad \quad f^2:=e^{-x}(u_{xx}\sin y-u_{xy}\cos y),\\
f^3&:=u-xu_x+xu_{xx}-u_{xx}+yu_{xy}-yu_{y},\quad f^4:=u_y-u_{xy},\quad f^5:=u_{xx}-u_x.
\end{align*}
This projects to
$$u=c_3+c_4y-c_5x+e^x\frac{c_2\sin y-c_1\cos y}{\sin^2 y-\cos^2 y}.$$

\end{example}
\subsection{The group invariance method}
In this section we discuss the group invariance approach to solve the differential conditions on the $\phi^\alpha$:
\begin{equation}\label{eq:diff cdtns}
 d\phi^\alpha\wedge\Omega_\theta\wedge\Omega_\phi=0
\end{equation}
We apply steps $2$ and $3$ of the algorithm together.

Firstly we explain step $3$ in this context. In applying Theorem \ref{sherr result}, we need five symmetries in a solvable structure. Let $\bar{\Omega}=\Omega_\theta\wedge\Omega_\phi$ be a $5-$form, where $\Omega_\phi$ comes from the algebraic conditions. We wish $\bar{\Omega}$ to be Frobenius integrable and ker $\bar{\Omega}=\Sp\{\bar{V_1},~\bar{V_2}\}$ -- the reduced Vessiot distribution. Now suppose there exists a solvable structure of $5$ linearly independent vector fields $X_1,\dots,X_5$. We will impose the condition $\bar\Omega(X_1,X_2,X_3,X_4,X_5)\neq 0$ and the conditions that $X_1$ is a symmetry of $\bar{\Omega},$ $X_2$ is a symmetry of $X_1\righthalfcup\bar{\Omega}$ and so on down to $X_5$ being a symmetry of $X_4\righthalfcup X_3\righthalfcup X_2\righthalfcup X_1\righthalfcup\bar{\Omega}.$

These symmetry conditions are
\begin{align}\label{eq:sym conds}
  \pounds_{X_1}\bar{\Omega}&=\lambda_1\bar{\Omega},\nonumber\\
  \pounds_{X_2}(X_1\righthalfcup\bar{\Omega})&=\lambda_2(X_1\righthalfcup\bar{\Omega}),\\
  \vdots\nonumber\\
\pounds_{X_5}(X_4\righthalfcup X_3\righthalfcup X_2\righthalfcup X_1\righthalfcup\bar{\Omega})&=\lambda_5(X_4\righthalfcup X_3\righthalfcup X_2\righthalfcup X_1\righthalfcup\bar{\Omega}),\nonumber
\end{align}
for some smooth functions $\lambda_1,\ldots,\lambda_5.$ Note that at this stage we have not assumed the Frobenius integrability of $\bar\Omega.$

The way in which the symmetry conditions act as a catalyst to the solution of the differential conditions on $\Omega_\phi$ can be seen more comprehensively by recalling that the sequence
\begin{equation}
\bar\Omega,\ X_1\righthalfcup\bar\Omega,\ X_2 \righthalfcup X_1\righthalfcup\bar\Omega,\ldots,
\ X_4\righthalfcup\ldots X_1\righthalfcup\bar\Omega
\end{equation}
are all simple forms and that they are all Frobenius integrable if $\bar\Omega$ is Frobenius integrable. Hence the Frobenius integrability of each form in the sequence is a necessary condition for the Frobenius integrability of $\bar\Omega$. And since each form in the sequence has one-form factors in $Sp\{\theta^a,\phi^\alpha\}$ their Frobenius integrability  represents a successive simplification of the Frobenius integrability of $\bar\Omega$.

We will now assume further that $X_1, X_2, X_3$ are linearly independent symmetries of  $\Omega_\theta$ satisfying $X_i\righthalfcup \Omega_\theta\neq 0,$ that is, symmetries of our PDE.  And we insist that the  first three symmetry conditions of \eqref{eq:sym conds} hold. For example, this means that
\begin{equation}
\pounds_{X_1}\Omega_\phi \equiv \mu_1\Omega_\phi \quad (\text{mod}\ \theta^a)
\end{equation}
which imposes conditions on the $\phi^\alpha$. This is the essence of the technique.

The remaining two symmetries $X_4, X_5$ in \eqref{eq:sym conds} of necessity have non-zero components, $X^V_4, X^V_5$, in the kernel of $\Omega_\theta.$ If $\text{ker}\ \bar\Omega \bigoplus \text{Sp}\{X_1,X_2,X_3,X_4\}$ is Frobenius integrable then the one-form
\begin{equation}\label{Eq:X4X5omega}
  \omega:=\frac{X_4\righthalfcup (X_3\righthalfcup X_2\righthalfcup X_1\righthalfcup (\Omega_\theta\wedge\Omega_\phi))}{(\Omega_\theta\wedge\Omega_\phi)(X_1,\ldots,X_5)},
\end{equation}
is closed by virtue of \eqref{eq:sym conds} and it is a constraint form for $\bar\Omega$. This closure is a necessary condition for the integrability of $\text{ker}\ \bar\Omega.$ We can simplify this differential condition by writing $\omega$ as
\begin{align*}
  \omega&=\frac{\Omega_\theta(X_1,X_2,X_3)X^V_4\righthalfcup\Omega_\phi}
  {\Omega_\theta(X_1,X_2,X_3)\Omega_\phi(X^V_4,X^V_5)}\\
  &=\frac{X^V_4\righthalfcup\Omega_\phi}{\Omega_\phi(X^V_4,X^V_5)} \in \Sp\{\phi^1,~\phi^2\}.
\end{align*}
(It should be noted that this {\em does not} imply that
$\frac{X^V_5\righthalfcup\Omega_\phi}{\Omega_\phi(X^V_4,X^V_5)}\quad \text{is closed} \mod~\omega).$\newline
In forcing the symmetry conditions \eqref{eq:sym conds} and the closure of $\omega$ we have found a closed form in $\Omega_\phi$ and reduced the two differential conditions \eqref{eq:diff cdtns} to a single condition
\begin{equation}
 d\phi\wedge\Omega_\theta\wedge\Omega_\phi=0.
\end{equation}
(Note that $\Omega_\phi=\phi^1\wedge\phi^2$.)
This result can clearly be generalised to any dimension but in general it isn't true that $Sp\{\phi^\alpha\}$ contains closed forms, even in this case. This is because generically we don't have a solvable structure. The least we can expect is given in the following proposition
\begin{proposition}
  Suppose that $\phi^{\alpha}$ satisfies
  \begin{align*}
    d\theta^a\wedge\Omega_\theta\wedge\Omega_\phi=0,\quad a=1,...,N,\\
    d\phi^{\alpha}\wedge\Omega_\theta\wedge\Omega_\phi=0, \quad \alpha=1,...,m.
  \end{align*}
 Then there exist $\tilde{\phi}^\alpha=l^\alpha_\beta\phi^\beta,~~|l^\alpha_\beta|\neq 0$ and $h^\alpha_a$ such that
 \begin{equation}
   df^\alpha=\tilde{\phi}^\alpha+h^\alpha_a\theta^a,\quad \alpha=1,...,m.
 \end{equation}
  \begin{proof}
Suppose that $\Sp\{\theta^a,\phi^\alpha\}$ is Frobenius integrable,
\begin{align*}
    d\theta^a\wedge\Omega_\theta\wedge\Omega_\phi=0,\quad a=1,...,N,\\
    d\phi^{\alpha}\wedge\Omega_\theta\wedge\Omega_\phi=0, \quad \alpha=1,...,m.
  \end{align*}
 This implies that \[df^a\in\Sp\{\theta^a,\theta^\alpha\},\quad \text{with}\quad df^1\wedge...\wedge df^A\neq 0.\]
 So, $df^A=h^A_a\theta^a+l^A_\alpha\phi^\alpha,\quad \text{for some} \quad h^A_a,~l^A_\alpha.$ By relabelling the $f^A$ and the $\phi^\alpha$ we have
  \begin{align*}
    df^1&=h^1_a\theta^a+l^1_\alpha\phi^\alpha,\\
    \vdots\\
     df^m&=h^m_a\theta^a+l^m_\alpha\phi^\alpha,\quad \text{with}\quad |l^\alpha_\beta|\neq0.~(\text{since}~\phi^\alpha\in\Sp\{df^A\})
  \end{align*}
 Then put $\tilde{\phi}^\alpha=l^\alpha_\beta\phi^\beta,\quad \alpha=1,\ldots,m.$ So that $df^\alpha=h^\alpha_a\theta^a+\tilde{\phi}^\alpha,\quad \alpha=1,\ldots,m.$
  \end{proof}
\end{proposition}We close this section with the following example.
\begin{example}
Consider the Gibbons-Tsarev equation \cite{Sergey} \begin{equation}\label{GT equation}
u_{yy}=u_xu_{xy}-\beta u_yu_{xx}+\mu u+\nu.
\end{equation}
We follow these authors by considering $\beta=1,~\mu=0$ and $\nu=1.$ Now apply the step $1$ of the algorithm.
The solution of the algebraic equation \eqref{eqn:FI1} implies
\begin{align*}
a_2=b_1,~b_2= - a_1u_y +b_1u_x.
\end{align*}
By substituting these values in the equation \eqref{eqn:phi1,2}, we have
  \begin{align*}
  \phi^1=\psi^3-a_1\psi^1-b_1\psi^2,\quad \quad \phi^2=\psi^4-b_1\psi^1-(-a_1u_y+b_1u_x)\psi^2.
\end{align*}
At this stage we have two unknown functions $a_1$ and $b_1.$ We impose the symmetry conditions \eqref{eq:sym conds} on $\bar{\Omega}.$ The symmetries of the original PDE \eqref{GT equation} are
\begin{align*}
X_1&:=\vp{x},\quad \quad X_2:=\vp{u},\quad \quad X_3:=\vp{y},\\
X_4&:=-2y\vp{y}-3x\vp{x}-4u\vp{u}-2u_y\vp{u_{y}}-u_x\vp{u_{x}}+2u_{xx}\vp{u_{xx}}+u_{xy}\vp{u_{xy}},\\
X_5&:=x\vp{y}-2x\vp{u}-2\vp{u_x}-u_x\vp{u_y}-u_{xx}\vp{u_{xy}}.
\end{align*}
The first three symmetry conditions \eqref{eq:sym conds} implies
\begin{align*}
  &\pd{a_1}{x}=0,~ \pd{a_1}{y}=0,~ \pd{a_1}{u}=0,~ \pd{b_1}{x}=0,~ \pd{b_1}{y}=0,~ \pd{b_1}{u}=0
\end{align*}
This means that $a_1,~b_1$ are independent of $x,~y,~u.$
The $4^{th}$ and $5^{th}$ symmetry conditions generate further $8$ conditions. We have $8$ first order partial differential equations for two unknown functions of $4$ variables. The solution of this overdetermined system is
\[a_1=0,~~b_1=0.\]
So, $\phi^1=\psi^3,~~\phi^2=\psi^4$ satisfy the differential conditions and the corresponding distribution is Frobenius integrable.  For this case we do not need to impose the further necessary condition on $\Omega_\phi.$

The corresponding reduced Vessiot distribution is
\begin{align*}
\bar{V}_1&:=\vp{x}+u_x\vp{u}+u_{xx}\vp{u_x}+u_{xy}\vp{u_y},\\
\bar{V}_2&:=\vp{y}+u_y\vp{u}+u_{xy}\vp{u_x}+(u_xu_{xy}-u_yu_{xx}+1)\vp{u_y}.
\end{align*}
After integration, we obtain
\begin{align*}
f^1&:=x+y\frac{u_{xy}}{u_{xx}}-\frac{u_{x}}{u_{xx}},\\
f^2&:=-u+xu_x-xyu_{xy}-\frac{1}{2}x^2u_{xx}+\frac{1}{2u_{xx}^3}(2yu_xu_{xy}u_{xx}^2-y^2u_{xy}^2u_{xx}^2-2yu_{xy}^2u_{xx}\\&
+2u_xu_{xy}u_{xx}^2-2u_{xy}^2-2u_{y}u_{xx}^2+2u_{xx}),\\
f^3&:=y+\frac{1}{u_{xx}}\log (u_x\frac{u_{xy}}{u_{xx}}-\frac{u_{xy}^2}{u_{xx}^2}+\frac{1}{u_{xx}}-u_y),\quad f^4:=-\frac{1}{2}\log u_{xx},\quad \quad f^5:=\frac{u_{xy}}{u_{xx}}.
\end{align*}
This projects to
$$u=-c_2+xyc_4c_5-xc_1c_4-yc_1c_4c_5-yc_5^2+c_4^{-1}(y+e^{c_3c_4-yc_4})+\frac{1}{2}c_4(x^2+y^2c_5^2).$$
\end{example}
\section{Summary}

In this paper, we have presented an algorithmic approach to the solution of second order partial differential equations of one dependent variable and two independent variables. We have successfully applied this algorithm to a wide class of linear and non-linear PDEs. Vessiot theory is used to generate the integrable sub-distributions  which satisfy the appropriate independence condition. Then a solvable structure of each integrable sub-distribution is found (by using DIMSYM for example). For those solvable structures of sufficient dimension, we integrate the reduced Vessiot distribution giving a local solution of the original system as the parameterised integral submanifold of the sub-distribution. These solutions appear in the original coordinates of the problem. Different sub-distributions correspond to distinct solutions. In this process we subdivide the Frobenius condition into two parts, one being algebraic and one differential. Solving the former is straight forward, the most challenging task is to satisfy the differential conditions (step $2$).  For the situations where the differential conditions are not satisfied easily, we have given two approaches to solve the problem. One is an application of the standard EDS method and other involves the simultaneous use of a solvable structure of symmetries to solve the conditions and to explicitly obtain the integrable submanifolds, giving a (generalised) group invariant solution. We give concrete examples throughout the paper.

Clearly future tasks involve the generalisation to higher order and higher dimensional systems, but we believe that the promise of symmetry reduction in the original coordinates and the calculus that is available for examining the shape of the solutions are exciting advances.

\section*{Acknowledgements}

Naghmana Tehseen gratefully acknowledges the support of a La Trobe University  postgraduate research award and the kind hospitality of Thomas Ivey. This work was supported by the Australian Research Council through a Discovery Project Grant $DP1095044$. For their helpful discussions the authors are grateful to Peter Vassiliou, Thomas Ivey, Phil Broadbridge and Dimetre Triadis.


\providecommand{\bysame}{\leavevmode\hbox to3em{\hrulefill}\thinspace}
\providecommand{\MR}{\relax\ifhmode\unskip\space\fi MR }
\providecommand{\MRhref}[2]{%
  \href{http://www.ams.org/mathscinet-getitem?mr=#1}{#2}
}
\providecommand{\href}[2]{#2}

\end{document}